\documentclass[12pt,reqno]{amsart}

\usepackage{amssymb}
\usepackage{amsmath}
\usepackage{amsthm,cite}

\newtheorem{lemma}{Lemma}[section]

\newtheorem{theorem}{Theorem}[section]
\newtheorem{corollary}{Corollary}[section]

\makeatletter \@addtoreset{equation}{section} \makeatother

\numberwithin{equation}{section}

\newcommand\lam{\lambda}
\newcommand\na{\nabla}
\newcommand\al{\alpha}

\newcommand\Del{\Delta}

\newcommand\f{\frac}
\newcommand\va{\varepsilon}

\newcommand\Ric{\textup{Ric}}

\begin{document}
\bibliographystyle{plain}
\title[]
{Gradient estimates for the Allen-Cahn equation  on Riemannian manifolds}

\author{Songbo Hou}
\address{Department of Applied Mathematics, College of Science, China Agricultural
University,  Beijing, 100083, P.R. China}
\email{housb10@163.com}

\date{}
\def\baselinestretch{1}

\begin{abstract}
In this paper, we consider bounded positive solutions to the Allen-Cahn equation on complete noncompact Riemannian manifolds without boundary. We derive gradient estimates for those solutions. As an application,  we get a Liouville type theorem on manifolds with nonnegative Ricci curvature.

\end{abstract}

\maketitle \baselineskip 18pt

\section{Introduction}
The  Allen-Cahn equation
\begin{align}\Del u+(1-u^2)u=0, \end{align}
   has its origin in the gradient theory
of phase transitions \cite{AC79}, and has attracted a lot of attentions in the last decades.   The famous De Giorgi conjecture  states that for $n\leq 8$, any entire solution to (1.1) in $\mathbb{R}^n$ with $|u|<1$ which is monotone in one
direction should be one-dimensional \cite{DE78}.  The conjecture was proved in dimension 2 by Ghoussoub-Gui \cite{GC} and
in dimension 3 by Ambrosio-Cabr\'{e} \cite{AC}, and in dimensions $4 \leq n \leq 8$ by Savin
\cite{SA}, under an extra
assumption. For $n\geq 9$, the conjecture is false \cite{DKW}.

Solutions to the Allen-Cahn equation have  the intricate
connection to the minimal surface theory. There are many results in the literature,  such as,  solutions concentrating along non-degenerate, minimal hypersurfaces
of a compact manifold were found in \cite{FR}.  So the equation is also an interesting topic for geometry.

The gradient estimate is a  useful method in the study of elliptic and parabolic equations.
It was originated by Yau \cite{Y75}, Cheng-Yau\cite{CY75}, and Li-Yau \cite{LY86}, and was extended by many authors, say Li\cite{Li91}, Negrin\cite{N95}, Souplet-Zhang \cite{SZ06}, Ma \cite{Ma}, Yang \cite{Yang,Yang10}, Cao \cite{Cao13} for various purposes.
   In this paper, we consider  bounded  positive solutions to Eq.(1.1)  and get the following theorem.

\begin{theorem} Let $M$ be a complete  noncompact $n$-dimensional   Riemannian manifold without boundary. Denote by  $B_p(2R)$ the geodesic ball of radius $2R$ around $P\in M$. Suppose $\Ric\geq -K(2R)$ in $B_p(2R)$  with $K(2R)\geq 0$,  $u$ is a bounded positive smooth solution of (1.1) on $M$ $u \leq C$ where $C$ is a positive constant.

(1) If $C\leq 1$, then we have

\begin{align*}
&\f{|\na u|^2}{u^2}+\frac{2}{3}(1-u^2)\\
&\leq\f{n}{1-\va}\left(\frac{ 2C_1^2+(n - 1)C_1^2(1+R\sqrt{K(2R)})+C_2}{{{R^2}}}+ \f{2n}{(1-\va)}\f{C_1^2}{R^2}+2K(2R)\right)
\end{align*}
on $B_p(R)$, where $C_1$, $C_2$ are positive constants, $0<\va<1$.

(2) If $C> 1$, then we have
\begin{align*}
&\f{|\na u|^2}{u^2}+s(1-u^2)\\
&\leq \f{ns^2}{2(1-\va)}\left(\frac{n}{{4\left( {1 - \varepsilon } \right)}}\frac{s^2}{{\left( {sq+s-1} \right)}}\frac{C_1^2}{R^2}
\frac{ 2C_1^2+(n - 1)C_1^2(1+R\sqrt{K(2R)})+C_2}{{{R^2}}}\right)\\
&+\f{ns^2}{2(1-\va)(s-1)}K(2R)+\f{s}{q}\sqrt{\f{n}{2(1-\va)}}C^2
\end{align*}
on $B_p(R)$, where $C_1$, $C_2$ are positive constants; $0<\va<1$, $s>1$, $q>0$ such that $\f{2(1-\va)}{n}\f{s-1}{sq}\ge \f{1}{\va}-1+\f{(3s-1)^2}{2}$.
 In particular, we can choose
$q=\f{2(1-\va)(s-1)}{ns\left[\f{1}{\va}-1+\f{(3s-1)^2}{2}\right]}$. Taking  $s=2$ and $\va=1/2$, we get
\begin{align*}
\f{|\na u|^2}{u^2}+2(1-u^2)&\leq 4n\left(\frac{54n^{2}}{27n+2}\frac{C_1^2}{R^2}
+\frac{ 2C_1^2+(n - 1)C_1^2(1+R\sqrt{K(2R)})+C_2}{{{R^2}}}\right)\\
&+4nK(2R)+54n\sqrt{n}C^2.
\end{align*}
 \end{theorem}

 As a consequence of Theorem 1.1, we have the following:

\begin{corollary}
Let $M$ be a complete noncompact $n$-dimensional  Riemannian manifold with Ricci tensor $\Ric\geq -{k }\,({k\geq 0})$.
 Suppose $u$ is a positive solution of (1.1) and $u\leq C$.

(1) If $C\leq 1$, we have
\begin{align*}
\f{|\na u|^2}{u^2}+\frac{2}{3}(1-u^2)\leq\f{2nk}{1-\va}.
\end{align*}
Letting $\va$ approach zero, we get
$$\f{|\na u|^2}{u^2}+\frac{2}{3}(1-u^2)\leq 2nk.$$
Furthermore,
$$|\na u|^2\leq 2nk.$$
(2) If $C>1$, we have
\begin{align*}
\f{|\na u|^2}{u^2}+s(1-u^2)\leq \f{ns^2k}{2(1-\va)(s-1)}+\f{s}{q}\sqrt{\f{n}{2(1-\va)}}C^2.
\end{align*}
In particular, choosing $s=2$ and $\va=1/2$, we have
\begin{align*}
\f{|\na u|^2}{u^2}\leq 4nk+(54n\sqrt{n}+2)C^2.
\end{align*}
Furthermore,
$$|\na u|^2 \leq (4nk+(54n\sqrt{n}+2)C^2)C^2.$$
\end{corollary}

For an application of Corollary 1.1, we get the following Liouville type theorem:
\begin{theorem}
Let $M$ be a complete  noncompact $n$-dimensional   Riemannian manifold with nonnegative Ricci curvature. If $u$ is a solution of (1.1)
with $0<u\leq 1$, then $u$ is equal to $1$ identically on $M$.
\end{theorem}

 In general, let $F\in C^2(\mathbb{R})$ be a nonnegative function and $u\in C^3({\mathbb{R}^n})$ a bounded entire solution  in  $\mathbb{R}^n$ of the equation
 $$\Del u=f(u),$$
 where $f=F^{'}$ is the first derivative of $F$. L. Modica \cite{LM85} proved that $|\na u|^2(x)\leq 2F(u(x))$ for every $x\in \mathbb{R}^{n}$.
 Later
Ratto-Rigoli \cite{RR95} extended Modica's result to  manifolds with nonnegative Ricci curvature. Also the conclusion of Theorem 1.2 can be
deduced from the result of Ratto and Rigoli by setting $F(u)=\frac{1}{4}u^4-\frac{1}{2}u^2+\frac{1}{4}$.
 However our result gives an explicit bound of $|\na u|$ in the case $k\neq 0$. In addition, Corollary 1.1 implies that the equation (1.1) does not admit an entire solution with values in $(0, 1)$ on manifolds with nonnegative Ricci curvature. The method in this paper can be also applied  to the equation
 $$\Del u+u^p-u^q=0,$$
 where $p,q\in \mathbb{R}$.

The rest of the paper is arranged as follows. In Section 2, we prove a basic lemma. In Section 3, we prove the main results.

\section{ Basic Lemma}

We consider
 $$W(x)=u^{-q}, $$ as the one defined in \cite{Li91},  where $q$ is a  positive constant to be chosen  later. A straightforward computation shows that

 $$\nabla W=-qu^{-q-1}\na u,$$
 $$|\nabla W|^2=q^2u^{-2q-2}|\na u|^2,$$
  \begin{equation}\frac{|\na W|^2}{W^2}=q^2u^{-2}|\na u|^2, \end{equation}

\begin{equation}\begin{aligned}\Del W&=q(q+1)u^{-q-2}|\na u|^2-qu^{-q-1}\Del u\\
 &=\frac{q+1}{q}\frac{|\na W|^2}{W}+qW-qW^{\f{q-2}{q}}.
 \end{aligned}
 \end{equation}

 We introduce the function
 \begin{equation}F(x)=\frac{|\na W|^2}{W^2}+\alpha (1-W^{-2/q}),\\
\end{equation}
 where $\al$ is a  positive constant to be fixed later.

Now we calculate
 \begin{equation}\na F(x)=\frac{\na |\na W|^2}{W^2}-\frac{2|\na W|^2\na W}{W^3}+\frac{2\alpha }{q}W^{-(q+2)/q}\na W,\end{equation}

 \begin{equation}
 \begin{aligned}
 \Delta F(x)&=\frac{2|\na^2 W|^2}{W^2}+\f{2\langle\na W,\Del \na W\rangle}{W^2}-8\frac{\langle\na^2 W, \na W \otimes \na W\rangle}{W^3}\\
 &+6\frac{|\na W|^4}{W^4}-2\frac{|\na W|^2\Del W}{W^3}
-\frac{2\alpha (q+2)}{q^2}\frac{W^{-2/q}}{W^2}|\na W|^2\\
&+\frac{2\alpha }{q}W^{-(q+2)/q}\Del W.
\end{aligned}
\end{equation}
Noting (2.2)  we have
 \begin{equation}\begin{aligned}\frac{2\langle\na W,\Del \na  W\rangle}{{{W^2}}} &= \frac{2\langle\na W, \na \Del W\rangle}{W^2} + \frac{2\Ric\langle\na W, \na W\rangle}{{{W^2}}}\\
&=\frac{4(q+1)}{q}\frac{\langle\na^2 W, \na W \otimes \na W\rangle}{W^3}-\frac{2(q+1)}{q}\frac{|\na W|^4}{W^4}\\
&+\frac{2|\na W|^2}{W^2}\left[q-(q-2)W^{-2/q}\right]+\frac{2\Ric\langle\na W, \na W\rangle}{{{W^2}}}, \end{aligned} \end{equation}

 \begin{equation}
 - \frac{2|\na W|^2\Del W}{W^3}
  {\rm{ =  - }}\frac{{2(q + 1)}}{q} \frac{|\na W|^4}{{{W^4}}} - 2q\frac{|\na W|^2}{W^2} + 2q{W^{ - \frac{2}{q}}}\frac{|\na W|^2}{W^2},
   \end{equation}
  \begin{equation}\frac{2\alpha }{q}W^{-(q+2)/q}\Del W=\frac{2\al (q+1)}{q^2}W^{-2/q}\frac{|\na W|^2}{W^2}+2\al W^{-2/q}-2\al W^{-4/q}.
    \end{equation}
By the H\"{o}lder inequality, we have
\[\frac{2\varepsilon |\na ^2 W|^2}{{{W^2}}} + \frac{2}{\varepsilon } \cdot \frac{|\na W|^4}{{{W^4}}} \ge 4\frac{\langle\na^2 W, \na W \otimes \na W\rangle}{{{W^3}}}.\]
Hence
\begin{align*}\frac{2|\na ^2 W|^2}{{{W^2}}} - 8\frac{\langle\na^2 W, \na W \otimes \na W\rangle}{{{W^3}}} + 6\frac{|\na W|^4}{{{W^4}}}& \ge \frac{{2\left( {1 - \varepsilon } \right)|\na^2 W|^2}}{{{W^2}}} - 4\frac{\langle\na^2 W, \na W \otimes \na W\rangle}{{{W^3}}}\\
 &+\left( {6 - \frac{2}{\varepsilon }} \right)\frac{|\na W|^4}{{{W^4}}},\end{align*}
where $0<\varepsilon<1.$

Using the fact  $ |\na ^2 W|^2 \ge \frac{1}{n}{\left( \Del W\right)^2}$,  we get
 \begin{equation}
 \begin{aligned}\frac{2|\na ^2 W|^2}{{{W^2}}} &- 8\frac{\langle\na^2 W, \na W \otimes \na W\rangle}{{{W^3}}} + 6\frac{|\na W|^4}{{{W^4}}}\\
  &\ge \frac{2\left(1 - \varepsilon  \right)}{n} \left(\f{\Del W}{W}\right)^2- 4\left(\frac{\langle\na^2 W, \na W \otimes \na W\rangle}{{{W^3}}}-\f{|\na W|^4}{W^4}\right) \\
  &-2\left( {\frac{1}{\varepsilon }-1} \right)\frac{|\na W|^4}{{{W^4}}}.
  \end{aligned} \end{equation}
 By (2.4),
 \begin{equation}\nabla {F} \cdot \nabla \log W = \frac{2\langle\na^2 W, \na W \otimes \na W\rangle}{{{W^3}}} - \frac{2|\na W|^4}{{{W^4}}} + \frac{{2\alpha }}{q}{W^{ - \frac{2}{q}}}\frac{|\na W|^2}{{{W^2}}}. \end{equation}
From (2.5) to (2.10),  we obtain
 \begin{equation}
 \begin{aligned}
 \Delta F  &\ge \frac{{2\left( {1 - \varepsilon } \right)}}{n}\left(\f{\Del W}{W}\right)^2 - 2\left( {\frac{1}{\varepsilon } - 1} \right)\frac{|\na W|^4}{{{W^4}}}\\
  &+ \f{2}{q}\langle\nabla {F}, \nabla \log W\rangle\\
  &+\left(4-6\f{\al }{q^2}\right)\frac{|\na W|^2}{{{W^2}}}{W^{ - \frac{2}{q}}}\\
   &+ \frac{2\Ric\langle\na W, \na W\rangle}{{{W^2}}}+2\al W^{-\f{2}{q}}(1-W^{-\f{2}{q}}).
\end{aligned}
 \end{equation}
It follows from  (2.2) and (2.3) that
 \begin{equation}\frac{{\Delta W}}{W} = \frac{q}{\alpha }F + \left( {\frac{{q + 1}}{q} - \frac{q}{\alpha }} \right)\frac{{{{\left| {\nabla W} \right|}^2}}}{{{W^2}}}. \end{equation}
 Set $\alpha=sq^2$, then
 \begin{equation}\frac{{\Delta W}}{W}{\rm{ = }}\frac{1}{{sq}}F + \left( {\frac{{q + 1}}{q} - \frac{1}{{sq}}} \right)\frac{{{{\left| {\nabla W} \right|}^2}}}{{{W^2}}} = \frac{1}{{sq}}F + {\left( {\frac{{q + 1 - 1/s}}{q}} \right)}\frac{{{{\left| {\nabla W} \right|}^2}}}{{{W^2}}}.  \end{equation}
Substituting (2.13) into (2.11) gives
\begin{align*}
 \Delta F &\ge \frac{{2(1 - \varepsilon )}}{n}\frac{1}{{{s^2q^2}}} {F^2} + \left[ {\frac{{2(1 - \varepsilon )}}{n}\frac{{{{\left( {sq + s - 1} \right)}^2}}}{{{s^2}{q^2}}} - 2\left(\frac{1}{\varepsilon } - 1\right)} \right]\frac{|\na W |^4}{{{W^4}}}\\
 &+ \frac{{4(1 - \varepsilon )}}{n}\frac{{\left( {sq + s - 1} \right)}}{{{s^2}{q^2}}}F\frac{|\na W|^2}{{{W^2}}} + \f{2}{q}   \langle\nabla F, \na \log W\rangle\\
     &+\left(4-6s\right)\frac{|\na W|^2}{{{W^2}}}{W^{ - \frac{2}{q}}}\\
   &+\frac{2\Ric\langle\na W, \na W\rangle}{{{W^2}}}+2sq^2W^{-\f{2}{q}}(1-W^{-\f{2}{q}}).
 \end{align*}

We get the following lemma.

\begin{lemma} Let $M$ be a complete noncompact $n$-dimensional  Riemannian manifold without boundary.  If $F$ is defined by (2.3) where $\al =sq^2$,  then we have
 \begin{equation}\begin{aligned}
\Delta F &\ge \frac{{2(1 - \varepsilon )}}{n}\frac{1}{{{s^2q^2}}} {F^2} \\
&+ \left[ {\frac{{2(1 - \varepsilon )}}{n}\frac{{{{\left( {sq + s - 1} \right)}^2}}}{{{s^2}{q^2}}}- 2\left(\frac{1}{\varepsilon } - 1\right)} \right]\frac{|\na W|^4}{{{W^4}}} \\
&+ \frac{{4(1 - \varepsilon )}}{n}\frac{{\left( {sq + s - 1} \right)}}{{{s^2}{q^2}}}F\frac{|\na W|^2}{{{W^2}}}\\
 &+\f{2}{q}  \langle\nabla F, \log W\rangle+ 2{W^{ - \frac{2}{q}}}F \\
 &+ \frac{2\Ric\langle\na W, \na W\rangle}{{{W^2}}}+\left(2-6s\right)\frac{|\na W|^2}{{{W^2}}}{W^{ - \frac{2}{q}}}.
\end{aligned} \end{equation}
\end{lemma}

\section{Proof of  Main Results}
\vskip 10 pt
\noindent
{\it Proof of Theorem 1.1.}
 Chose a cut-off function  $\chi$  $ \in C^2[0,+\infty)$ such that $\chi (r)=1$ for $r\le 1$, $\chi (r)=0$ for $r>2$, and
$0\leq \chi (r)\leq 1$. In addition, we  require $\chi$ satisfies
$-C_1\leq \chi^{-1/2}(r)\chi'(r)\leq 0$  and $\chi''(r)\geq -C_2$,   where $C_1, C_2$ are positive constants.

For a fixed point $p$, denote by  $r(x)$  the geodesic distance between $x$ and $P$.  Define
\[\phi(x)=\chi\left(\f{r(x)}{R}\right).\]
 It is clear that  \[|\na \phi|^2\leq \f{C_1^2}{R^2}\phi.\]
By the Laplacian  comparison theorem,  we get\[\Del\phi \geq - \frac{ (n - 1)C_1^2(1+R\sqrt{K(2R)})+C_2}{{{R^2}}}. \]

  Now we consider the function $\phi(x)F(x)$. By the  argument of Calabi\cite{Ca58}, we assume that the function $\phi(x)F(x)$  is smooth in $B_P(2R)$. Let $z$ be the point where $\phi F$ achieves its maximum in $B_P(2R)$.   We can  assume that $\lam:=\phi(z)F(z)>0$ since the theorem is obviously true if $\lam \leq 0$.   Then we have
 \begin{equation}\nabla \left( \phi F \right)=\na \phi F+\phi \na F = 0 \end{equation}
and  \begin{equation}\Del (\phi F)\leq 0 \end{equation}
at the point $z$,

Using Eq.(3.1),  we have
\[\nabla F =  - \frac{{\nabla \phi}}{\phi}F.\]
By (3.2), we have
\[\Delta \phi \cdot F + 2\nabla \phi \cdot \nabla F + \phi\Delta F \le 0.\]
Thus we obtain
\[F\Delta \phi + \phi \Del F-2F\phi^{-1}|\na \phi|^2 \le 0\] at $z$.

Then  for
\[B=\frac{ 2C_1^2+(n - 1)C_1^2(1+R\sqrt{K(2R)})+C_2}{{{R^2}}},\] we have
\[\phi\Del F\leq BF.\]
Multiplying both sides of (2.14) by $\phi^2$, we obtain at $z$,
 \begin{equation}
\begin{aligned}
B\phi F &\ge \frac{{2(1 - \varepsilon )}}{n}\frac{1}{{{s^2q^2}}} {(\phi F)^2} \\
&+ \phi^2\left[ {\frac{{2(1 - \varepsilon )}}{n}\frac{{{{\left( {sq + s - 1} \right)}^2}}}{{{s^2}{q^2}}}- 2\left(\frac{1}{\varepsilon } - 1\right)} \right]\frac{|\na W|^4}{{{W^4}}} \\
&+ \phi ^2\frac{{4(1 - \varepsilon )}}{n}\frac{{\left( {sq + s - 1} \right)}}{{{s^2}{q^2}}}F\frac{|\na W|^2}{{{W^2}}}\\
 &+\f{2}{q}\phi^2  \langle\nabla F, \log W\rangle+ 2{W^{ - \frac{2}{q}}}\phi^2F \\
 &+ \frac{2\Ric\langle\na W,\na W\rangle}{{{W^2}}}\phi^2-\left(6s-2\right)\phi^2\frac{|\na W|^2}{{{W^2}}}{W^{ - \frac{2}{q}}}.
\end{aligned}
 \end{equation}
We consider two cases: (1) $C\leq 1$ and (2) $C>1$.

(1)
Since $u\leq 1$, it is easy to see that
 \begin{equation}\frac{2\Ric\langle\na W,\na W\rangle}{W^2}\phi ^2\geq -2K(2R)\f{|\na W|^2}{W^2}\phi^2\geq-2K(2R)\phi F \end{equation}
and
  \begin{equation}2{W^{ - \frac{2}{q}}}\phi^2F-\left(6s-2\right)\phi^2\frac{|\na W|^2}{{{W^2}}}{W^{ - \frac{2}{q}}}\geq -(6s-4)\phi F \end{equation}
  if $s\geq \frac{2}{3}$.

Substituting (3.4), (3.5) into (3.3), and  choosing $s=\f{2}{3}$ and $q>0$ small enough such that  $\f{(1-\va)}{n}\f{(sq+s-1)^2}{s^2q^2}\ge \f{1}{\va}-1$, then we have
 \begin{equation}
\begin{aligned}
B\phi F &\ge \frac{{9(1 - \varepsilon )}}{2nq^2}{(\phi F)^2} \\
&-\frac{{3(1 - \varepsilon )}}{nq^2}(\phi F)^2\\
 &-\f{2}{q}\phi F  \langle\nabla \phi, \f{\na W}{W}\rangle-2K(2R)\phi F.
\end{aligned}
 \end{equation}
 We take the similar technique as in \cite{Ma}. Clearly,
 \begin{equation}-\frac{2}{q}F\phi \langle\nabla \phi,\frac{{\nabla W}}{W}\rangle \ge -\f{2C_1}{qR}(\phi F)^{3/2}.  \end{equation}
Combining (3.6) and (3.7), we arrive at

\begin{align*}
B\phi F &\ge \frac{{3(1 - \varepsilon )}}{2nq^2} {(\phi F)^2}-\f{2C_1}{qR}(\phi F)^{3/2}-2K(2R)\phi F.
\end{align*}
It follows that
\begin{align*}
B+ \f{2C_1}{qR}(\phi F)^{1/2}+2K(2R)&\ge \frac{{3(1 - \varepsilon )}}{2nq^2} {(\phi F)}.
\end{align*}
In other words, we get
 \begin{equation}
B+ \f{2C_1}{qR}\lam^{1/2}+2K(2R)\ge \frac{{3(1 - \varepsilon )}}{2nq^2} \lam.
 \end{equation}
Note that
 \begin{equation}\f{2C_1}{qR}\lam^{1/2}\leq \f{(1-\va)}{2nq^2}\lam+\f{2n}{(1-\va)}\f{C_1^2}{R^2}.  \end{equation}
Substituting (3.9) into (3.8), we get
\begin{align*}
B+ \f{2n}{(1-\va)}\f{C_1^2}{R^2}+2K(2R)&\ge \frac{{1 - \varepsilon }}{nq^2} \lam.
\end{align*}
Then we get
\begin{equation}
\begin{aligned}
\lam\leq &\f{nq^2}{1-\va}\left(B+ \f{2n}{(1-\va)}\f{C_1^2}{R^2}+2K(2R)\right)\\
=&\f{nq^2}{1-\va}\left(\frac{ 2C_1^2+(n - 1)C_1^2(1+R\sqrt{K(2R)})+C_2}{{{R^2}}}\right.\\
+ &\left.\f{2n}{(1-\va)}\f{C_1^2}{R^2}+2K(2R)\right).
\end{aligned}
\end{equation}
(2)  By the condition on Ricci curvature, we derive
\[\frac{2\Ric\langle\na W,\na W\rangle}{W^2}\phi ^2\geq -2K(2R)\f{|\na W|^2}{W^2}\phi^2.\]
By H\"{o}lder's inequality, we get
 \[2K(2R)\phi^2\frac{{{{\left| {\nabla W} \right|}^2}}}{{{W^2}}} \le \frac{{2\left( {1 - \varepsilon } \right)}}{n}\f{(s-1)^2}{s^2q^2}\frac{{{{\left| {\nabla W} \right|}^4}}}{{{W^4}}}\phi^2 + \frac{n}{2\left( {1 - \varepsilon } \right)}\f{s^2q^2}{(s-1)^2}{K^2}(2R)\phi^2\]
and
 \[(6s-2)\phi ^2\f{|\na W|^2}{W^2}W^{-2/q}\leq \f{(6s-2)^2}{4}\f{|\na W|^4}{W^4}\phi^2+C^4\phi^2.\]
By (3.4),
\[\f{2}{q}\phi^2 \langle \nabla F, \log W\rangle=-\frac{2}{q}\phi F \langle\nabla \phi,\frac{{\nabla W}}{W}\rangle.\]
Using H\"{o}lder's inequality again gives
\[\frac{2}{q} \phi F  \langle\nabla \phi,\frac{{\nabla W}}{W}\rangle \le \frac{{4\left( {1 - \varepsilon } \right)}}{n}\frac{{\left( {sq+s-1} \right)}}{{{s^2q^2}}}\frac{{{{\left| {\nabla W} \right|}^2}}}{{{W^2}}}F\phi^2 + \frac{n}{{4\left( {1 - \varepsilon } \right)}}\frac{s^2\phi F}{{\left( {sq+s-1} \right)}}\frac{{{{\left| {\nabla \phi} \right|}^2}}}{\phi}.\]
 Choose $s>1$ and $q>0$ such that  $\f{2(1-\va)}{n}\f{s-1}{sq}\ge \f{1}{\va}-1+\f{(3s-1)^2}{2}$.  Then (3.3) becomes
\begin{align*}
B\phi F &\ge \frac{{2(1 - \varepsilon )}}{n}\frac{1}{{{s^2q^2}}} {(\phi F)^2}
-\frac{n}{{4\left( {1 - \varepsilon } \right)}}\frac{s^2}{{\left( {sq+s-1} \right)}}\frac{C_1^2}{R^2}\phi F\\
&-\frac{n}{2\left( {1 - \varepsilon } \right)}\f{s^2q^2}{(s-1)^2}{K^2}(2R)-
C^4,\end{align*}
whence
\begin{align*}
0 &\ge \frac{{2(1 - \varepsilon )}}{n}\frac{1}{{{s^2q^2}}} {\lam^2}
-\left(\frac{n}{{4\left( {1 - \varepsilon } \right)}}\frac{s^2}{{\left( {sq+s-1} \right)}}\frac{C_1^2}{R^2}+B\right)\lam\\
&-\frac{n}{2\left( {1 - \varepsilon } \right)}\f{s^2q^2}{(s-1)^2}{K^2}(2R)-
C^4.\end{align*}
Thus
\begin{equation}
\begin{aligned}
\lam &\leq \f{ns^2q^2}{2(1-\va)}\left(\frac{n}{{4\left( {1 - \varepsilon } \right)}}\frac{s^2}{{\left( {sq+s-1} \right)}}\frac{C_1^2}{R^2}\right.
\\
+&\left.\frac{ 2C_1^2+(n - 1)C_1^2(1+R\sqrt{K(2R)})+C_2}{{{R^2}}}\right)\\
&+\f{ns^2q^2}{2(1-\va)(s-1)}K(2R)+sq\sqrt{\f{n}{2(1-\va)}}C^2.
\end{aligned}
\end{equation}

Combining (3.10) and (3.11), we conclude the theorem. $\hfill\Box$\\

{\it Proof of Corollary 1.1.} Passing to the limit $R\rightarrow +\infty$ in the estimates of Theorem 1.1,
we get the desired results. $\hfill\Box$\\

{\it Proof of Theorem 1.2.} Suppose that $M$ is a complete noncompact  Riemannian manifold with nonnegative Ricci curvature.
 If $u$ is a solution of (1.1) on $M$ and $0<u\leq 1$, then by
 Corollary 1.1, we get
\[\f{|\na u|^2}{u^2}+\frac{2}{3}(1-u^2)\leq 0.\]
It follows that $|\na u|\equiv 0$ and $u\equiv 1$. This concludes Theorem 1.2. $\hfill\Box$


\vskip 10pt






\vskip 30 pt
\begin{thebibliography}{10}

\bibitem{AC79} S. M. Allen, J. W. Cahn, A microscopic theory for antiphase boundary motion and its application to antiphase domain coarsening, Acta Metall. 27 (1979) 1085-1095.
\bibitem {AC} L. Ambrosio, X. Cabr\'{e}, Entire solutions of semilinear elliptic equations in $\mathbb{R}^3$ and a conjecture of
De Giorgi, J. Amer. Math. Soc. 13(4) (2000) 725-739 (electronic).
\bibitem{Ca58} E. Calabi, An extension of E. Hopf's maximum principle with application to Riemannian  geometry, Duke Math. J. 25 (1958) 45-46.
\bibitem{Cao13} X. Cao, B. Fayyazuddin Ljungberg, B. Liu, Differential Harnack estimates for a nonlinear heat equation, J. Funct. Anal. 265 (2013) 312-2330.
\bibitem{CY75}S. Y. Cheng, S. T. Yau, Differential equations on Riemannian manifolds and their geometric applications, Comm.
Pure Appl. Math. 28(3) (1975) 333-354.
\bibitem{DE78} E. De Giorgi, Convergence problems for functionals and operators,  In Proceedings of the International
Meeting on Recent Methods in Nonlinear Analysis (Rome, 1978), pages 131-188. Pitagora, Bologna,
1979.
\bibitem{DKW} M. del Pino, M. Kowalczyk, J. Wei, On De Giorgi's conjecture in dimension  $N\geq 9$, Ann. of Math. (2),
174(3) (2011) 1485-1569.
\bibitem{GC} N. Ghoussoub,   C. Gui,  On a conjecture of De Giorgi and some related problems,  Math. Ann.,
311(3) (1998) 481-491.
\bibitem {Li91} Jiayu Li, Gradient estimates and Harnack inequalities for nonlinear parabolic and nonlinear elliptic equations on Riemannian manifolds, J. Funct. Anal. 100 (1991) 233-256.
 \bibitem {Ma} L. Ma, Gradient estimates for a simple elliptic equation on non-compact Riemannian manifolds,
J. Funct. Anal.  241 (2006)  374-382.
\bibitem{LY86} P. Li, S. T. Yau, On the parabolic kernel of the Schr\"{o}dinger operator, Acta Math. 156 (1986) 153-201.

\bibitem{LM85} L. Modica, A gradient bound and a liouville theorem for nonlinear poisson equations, Commun. Pure Appl. Math. 38 (1985) 679-684.

\bibitem{N95} E. Negrin, Gradient estimates and a   Liouville type theorem for the Schr\"{o}dinger operator, J. Funct. Anal. 127 (1995) 198-203.
\bibitem{FR}  F. Pacard,  M. Ritor\'{e}, From the constant mean curvature hypersurfaces
to the gradient theory of phase transitions, J. Differential Geom. 64(3)
(2003) 356-423.

\bibitem{RR95} A. Ratto, M. Rigoli, Gradient bounds and Liouville's type theorems for the Poisson equation on complete Riemannian manifolds,
Tohoku Mathematical Journa l47(4) (1995) 509-519.
\bibitem {SA} O. Savin. Regularity of flat level sets in phase transitions, Ann. of Math.
(2) 169 (1) (2009) 41-78.
\bibitem{SZ06} P. Souplet, Q. S. Zhang, Sharp gradient estimate and Yau's Liouville theorem for the heat equation on noncompact manifolds. Bull. London Math. Soc. 38 (2006) 1045-1053.
\bibitem{Yang} Y. Yang,  Gradient estimates for a nonlinear parabolic equation on Riemannian
manifolds, Proc. Amer. Math. Soc.  136  (2008)  4095-4102.
\bibitem{Yang10} Y. Yang, Gradient estimates for the equation $\Delta u+cu^{-\alpha}=0$ on Riemannian manifolds, Acta Math. Sin. (Engl. Ser.)
26 (2010) 1177-1182.
\bibitem{Y75} S. T. Yau, Harmonic functions on complete Riemannian manifolds, Comm. Pure Appl. Math. 28 (1975) 201-228.



\end{thebibliography}
\end{document}